\documentclass[12pt,leqno]{article}
\usepackage{amsmath, amsthm, amsfonts, amssymb, color}
\usepackage{mathrsfs}
\theoremstyle{definition}

\newcommand{\scr}[1]{\mathscr #1}
\definecolor{wco}{rgb}{0.5,0.2,0.3}

\numberwithin{equation}{section} \theoremstyle{remark}

\newcommand{\ua}{\uparrow}

\title{{\bf Log-Harnack Inequality for   Gruschin Type Semigroups}\footnote{Supported in
 part by NNSFC(11131003), SRFDP, the Laboratory of Mathematical and  Complex Systems and the Fundamental Research Funds for the Central Universities.}
}
\author{
{\bf Feng-Yu Wang$^{1),2)}$ and  Lihu Xu$^{3)}$}\\
\footnotesize{$^{1)}$School of Mathematical Sciences,
Beijing Normal University, Beijing 100875, China}\\
  \footnotesize{$^{2)}$Department of Mathematics,
Swansea University, Singleton Park, SA2 8PP, UK}\\
\footnotesize{Email: wangfy@bnu.edu.cn; F.Y.Wang@swansea.ac.uk} \\
\footnotesize{$^{3)}$Department of Mathematics,
Brunel University, Uxbridge, UB8 3PH, UK}\\
\footnotesize{Email: Lihu.Xu@brunel.ac.uk}
}
\begin{document}
\def\R{\mathbb R}  \def\ff{\frac} \def\ss{\sqrt} \def\B{\mathbf
B}
\def\N{\mathbb N} \def\kk{\kappa} \def\m{{\bf m}}
\def\dd{\delta} \def\DD{\Delta} \def\vv{\varepsilon} \def\rr{\rho}
\def\<{\langle} \def\>{\rangle} \def\GG{\Gamma} \def\gg{\gamma}
  \def\nn{\nabla} \def\pp{\partial} \def\EE{\scr E}
\def\d{\text{\rm{d}}} \def\bb{\beta} \def\aa{\alpha} \def\D{\scr D}
  \def\si{\sigma} \def\ess{\text{\rm{ess}}}
\def\beg{\begin} \def\beq{\begin{equation}}  \def\F{\scr F}
\def\Ric{\text{\rm{Ric}}} \def\Hess{\text{\rm{Hess}}}
\def\e{\text{\rm{e}}} \def\ua{\underline a} \def\OO{\Omega}  \def\oo{\omega}
 \def\tt{\tilde} \def\Ric{\text{\rm{Ric}}}
\def\cut{\text{\rm{cut}}} \def\P{\mathbb P}
\def\C{\scr C}     \def\E{\mathbb E}
\def\Z{\mathbb Z}
  \def\Q{\mathbb Q}  \def\LL{\Lambda}
  \def\B{\scr B}    \def\ll{\lambda}
\def\vp{\varphi}

\maketitle
\begin{abstract} By constructing  a coupling  in two steps  and  using the Girsanov theorem under a regular conditional probability, the log-Harnack inequality is established for a large class of Gruschin type semigroups whose generator might be both degenerate and non-Lipschitzian.
\end{abstract} \noindent

 \vskip 2cm

\section{Introduction}

In recent years, regularity estimates has been investigated for some typical   subelliptic diffusion semigroups, see  \cite{GW, WZ, Z} for the study of generalized stochastic Hamiltonian systems,  and see \cite{BBBC, BGM, DM, Li} for gradient estimates and Harnack inequalities on  Heisenberg groups. This paper aims to investigate the log-Harnack   inequality introduced in \cite{RW09, W10} for   Gruschin type semigroups whose generators are degenerate and possibly singular.  This inequality is a weaker version of the dimension-free Harnack inequality introduced in \cite{W97}, and have a number of applications to   heat kernel estimates and transportation-cost inequalities, see e.g. \cite[Section 4]{W12}.

Let us start with the classical Gruschin semigroup on $\R^2$ with order $l>0,$ which is generated by
$$L(x^{(1)},  x^{(2)}):= \ff 1 2 \Big(\ff{\pp^2}{\pp (x^{(1)})^2} + |x^{(1)}|^{2l} \ff{\pp^2}{\pp (x^{(2)})^2}\Big).$$   The corresponding diffusion process can be constructed by solving the SDE
$$\beg{cases} \d X_t^{(1)}= \d B_t^{(1)},\\
\d   X_t^{(2)}= |X_t^{(1)}|^l\,\d  B_t^{(2)},\end{cases}$$ where $B_t:=(B_t^{(1)},  B_t^{(2)})$ is a two-dimensional Brownian motion.  Clearly, the equation is degenerate, and when $l<1$ the coefficient in the second equation is non-Lipschitzian.   In the simplest case that $l=1$,
the generalized curvature-dimension condition introduced  in \cite{BB3} holds, so that  the  gradient estimates and Harnack inequalities  derived in \cite{BB1,BB3} are valid for the associated semigroup. When $l$ is a natural number larger than $1$,  a more general version of curvature condition has been confirmed in \cite{W12}, which also implies  explicit gradient estimates of the semigroup.  Moreover, for general $l\ge 1$,  a Bismut type derivative formula was derived for the semigroup in \cite{W12a} by using Malliavin calculus.  However, due to the singularity of the coefficient,  the arguments used in these papers are no longer valid  if $l\in (0,1)$, and except for $l=1$, the log-Harnack inequality is not yet known for the semigroup.
In this paper we aim to establish the log-Harnack  inequality of the Gruschin semigroup for all $l>0.$ But, our argument used in the paper does not imply the dimension-free Harnack inequality in the sense of \cite{W97} for the Gruschin semigroup.

A key tool in the study is  the coupling method by change of measure  introduced in  \cite{ATW06}.
This method has been developed and applied to various finite- and infinite-dimensional models, see e.g.  \cite{GW,L,LW, O,ORW,Wan3,WY10,WWX10,WX10,X}   and references therein. Due to the high degeneracy (for large $l$) and the singularity (for small $l$) of the coefficient, we have to overcome new difficulty in the study.

We   consider the following   more general SDE  for $X_t:=(X_t^{(1)}, X_t^{(2)})$ on $\R^{m}\times \R^d= \R^{m+d} (m,d\ge 1)$:
\beq\label{EX} \beg{cases} \d X_t^{(1)}= b^{(1)}(t, X_t^{(1)}) \d t + \si^{(1)}(t)\, \d B_t^{(1)},\\
\d X_t^{(2)}= b^{(2)} (t,X_t)\d t + \si^{(2)}(t,X_t^{(1)})\,\d  B_t^{(2)},\end{cases}\end{equation}where $B_t:=(B_t^{(1)}, B^{(2)}_t)$ is the $(m+d)$-dimensional Brownian motion on a complete probability space $(\OO,\F,\P)$ with natural filtration $\{\F_t\}_{t\ge 0},$ and
\beg{equation*}\beg{split} & b^{(1)}: [0,\infty)\times\R^m\to \R^m,\ \ b^{(2)}: [0,\infty)\times \R^{m+d}\to \R^d,\\
 & \si^{(1)}: [0,\infty)\to \R^m\otimes\R^m,\ \ \si^{(2)}: [0,\infty)\times \R^m\to \R^d\otimes\R^d\end{split}\end{equation*} are measurable, and $b^{(1)}, b^{(2)}, \si^{(2)}$ are continuous in the second variable.
 Assume
 \beg{enumerate}
\item[{\bf (A1)}]  There exists a decreasing function $\ll: [0,\infty)\to (0,\infty)$ such that $\si^{(1)}(t)\si^{(1)}(t)^*\ge \ll_t^2 I_{m\times m},\ t\ge 0.$
\item[{\bf (A2)}] There exists an increasing function $K: [0,\infty)\to \R$ such that
$$\<b^{(1)}(t, x^{(1)})- b^{(1)}(t, y^{(1)}), x^{(1)}-y^{(1)}\> \le K_t |x^{(1)}-y^{(1)}|^2,\ t\ge 0, x^{(1)}, y^{(1)}\in \R^m.$$
\item[{\bf (A3)}]  There exist increasing functions $\Theta: [0,\infty)\to \R, h: [0,\infty)\to [1,\infty)$ and $\vp_\cdot: [0,\infty)^2\to [0,\infty)$ with $\vp(0)=0$ such that
\beg{equation*}\beg{split} & \<b^{(2)}(t, x)- b^{(2)}(t, y), x^{(2)}-y^{(2)}\> + \ff 1 2 \|\si^{(2)}(t,x^{(1)}) -\si^{(2)}(t, y^{(1)})\|_{HS}^2\\
& \le \Theta_t   |x^{(2)}-y^{(2)}|^2
 + \vp_t(|x^{(1)}-y^{(1)}|^2)h(|x^{(1)}|\lor |y^{(1)}|)\end{split}\end{equation*} holds for all $ t\ge 0$ and $ x=(x^{(1)}, x^{(2)}), y=(y^{(1)}, y^{(2)})\in \R^{m+d}.$
    \end{enumerate}

 It is well known that {\bf (A2)} implies the existence, uniqueness and non-explosion of strong solutions to the first equation in (\ref{EX}). Once
 $X_t^{(1)}$ is fixed,  then it follows from {\bf (A3)} that the second equation in (\ref{EX}) admits a unique global solution. Note that {\bf (A3)}
allows $\si^{(2)}(t,\cdot)$ to be merely H\"older continuous  when e.g. $\vp_t(r)= r^\aa$ for some constant $\aa\in (0,1)$.
 For any $x=(x^{(1)}, x^{(2)})\in \R^{m+d}$, we let $X_t(x)=(X^{(1)}_t(x), X^{(2)}_t(x))$ denote the solution to (\ref{EX}) with $X_0=x.$ Since $X^{(1)}_t(x)$ does not depend  on $x^{(2)}$ we also write $X^{(1)}_t(x)= X_t^{(1)}(x^{(1)}). $
 We intend to establish Harnack type inequalities for the associated semigroup $P_t$:
 $$P_tf(x):= \E f(X_t(x)),\ \ f\in \B_b(\R^{m+d}), t\ge 0, x\in\R^{m+d}.$$

We remark that {\bf (A1)} means that the first component process  $X^{(1)}_t$  is a non-degenerate diffusion process on $\R^m$, {\bf (A2)} is the usual semi-Lipschitz condition for this process, and when e.g. $b^{(2)}$ is independent of $x^{(1)}$ and semi-Lipschitzian in $x^{(2)}$, {\bf (A3)} holds provided
$$\ff 1 2 \|\si^{(2)}(t,x^{(1)}) -\si^{(2)}(t, y^{(1)})\|_{HS}^2\le \vp_t(|x^{(1)}-y^{(1)}|^2)h(|x^{(1)}|\lor |y^{(1)}|).$$ In particular, for the Gruschin semigroup where $\si^{(2)} (t, x^{(1)}) =|x^{(1)}|^l$,    this condition holds for $\vp_t(r)= r^{l\land 1}$ and $h(r)= c\lor r^{(l-1)^+}$ for some constant $c\ge 1.$

  In order to control the   degeneracy of $\si^{(2)}(t,\cdot)$, we need the condition
  \beq\label{PS} \psi_T (x^{(1)}, y^{(1)}):= \sup_{t\in [T, 2T]}  \E^{y^{(1)}} \Big\{ \big\|\si^{(2)}(t,X_t^{(1)})^{-1}\big\|^2\sup_{s\in [0,T]} h \big(|X_s^{(1)}|+|x^{(1)}-y^{(1)}|\big)\Big\}<\infty
  \end{equation} for  $T>0$ and $x^{(1)}, y^{(1)}\in \R^m,$ where $\E^{y^{(1)}}$ is the expectation for $X_t^{(1)}(y^{(1)}),$
   $\|\si^{-1}\|$ stands for the operator norm of the inverse of    a $d\times d$-matrix $\si$, and when the matrix is non-invertible we take $\|\si^{-1}\|=\infty.$

\beg{thm}\label{T1.1}  Assume that {\bf (A1)}, {\bf (A2)}, {\bf (A3)} and   $(\ref{PS})$ hold. Then for any strictly positive function $f\in\B_b(\R^{m+d}), x=(x^{(1)}, x^{(2)}), y=(y^{(1)}, y^{(2)})\in \R^{m+d}$ and $T>0$,
\beg{equation*}\beg{split}    P_{2T}\log f(y)\le &\log P_{2T} f(x)  +
   \ff{K_T|x^{(1)}-y^{(1)}|^2}{\ll_T^2(1-\e^{-2K_TT})}\\
   &+  \ff{\Theta_{2T}\e^{2\Theta_TT}\psi_T(x^{(1)}, y^{(1)})}{\e^{-2\Theta_{2T}T}-\e^{-4\Theta_{2T}T}}
  \bigg\{   |x^{(2)}-y^{(2)}|^2 +\ff{1-\e^{-2\Theta_TT}}{\Theta_T} \vp_T(|x^{(1)}-y^{(1)}|^2)\bigg\}.\end{split}\end{equation*}
\end{thm}

\

Let us come back to the classical Gruschin semigroup for which $m=d=1, b^{(1)}= b^{(2)}=0, \si^{(1)}= 1$ and $\si^{(2)}(t,x^{(1)})= |x^{(1)}|^l.$
Then {\bf (A1)}-{\bf (A3)} hold for $\ll=1, K= \Theta =0, \varphi (r)=r^{l\land 1}$ and $h(r)= c_1 \lor r^{(l-1)^+}$ for some constant $c_1\ge 1.$
When $l \in (0,\ff 1 2),$ we may take $h\equiv 1$ so that
$$ \psi_T(x^{(1)},y^{(1)})= \sup_{t\in [T,2T]}  \int_{\R } \ff 1 {|z|^{2 l}\ss{2t\pi} }\e^{-\ff{|z-x^{(1)}|^2}{2 t}}\d z
\le   \ff {c_2} {T^{l }}<\infty$$ for some constant $c_2>0$.   Therefore,   according to Theorem \ref{T1.1},  the log-Harnack inequality
$$P_{2T}\log f(y)\le \log P_{2T}f(x) + \ff{|x^{(1)}-y^{(1)}|^2}{2T} + \ff{c}{T^{l+1}}
  \Big\{|x^{(2)}-y^{(2)}|^2+2T|x^{(1)}- y^{(1)}|^{2l}\Big\}$$ holds. On the other hand, however, it is easy to see that $\psi_T=\infty$ for $l\ge \ff 1 2.$  Similarly, for the  Gruschin semigroup    on $\R^{m+d}$, i.e. $b^{(1)}=0, b^{(2)}=0, \si^{(1)}=I_{m\times m}$ and $\si^{(2)}(x^{(1)})= |x^{(1)}|^lI_{d\times d}$,   $\psi_T(x^{(1)}, y^{(1)})<\infty$ (and hence   log-Harnack inequality holds)   if and only if   $l\in (0,\ff m 2).$

\

To derive the log-Harnack inequality for the Gruschin semigroup for all $l>0$,  we intend to relax the condition (\ref{PS}) by using the invertibility of the following  integral matrix $Q_T$ to replace   that of $\si^{(2)}$. To this end,  we will need to assume that $b^{(2)}(t,x)$ is linear in $x^{(2)};$ that is, $b^{(2)}(t, x)= A x^{(2)}+ \tt b^{(2)}(t,x^{(1)})$ for some $d\times d$-matrix
$A$   and  some $\tt b^{(2)}\in C([0,\infty)\times \R^m;\R^d).$
  Let
$$Q_T= \int_T^{2T} \e^{A(T-t)} \si^{(2)} (t,X^{(1)}_t) \si^{(2)} (t,X^{(1)}_t)^*\e^{A^*(T- t)} \d t,\ \ T>0.$$

\beg{thm}\label{T1.2} Assume that  {\bf (A1)}, {\bf (A2)} and {\bf (A3)} hold for $b^{(2)}(t, x)= A x^{(2)}+ \tt b^{(2)}(t,x^{(1)}),$ where
$A$ is a $d\times d$-matrix and $\tt b^{(2)}\in C([0,\infty)\times \R^m;\R^d).$ Let $\theta_T= \sup_{t\in [0,T]}\|\e^{-At}\|.$ If $Q_T$ is invertible and
$$\Psi_T(x^{(1)},y^{(1)}):= \E^{y^{(1)}}  \bigg\{\big\|Q_T^{-1}\big\|^{2} \bigg(\int_T^{2T} \big\|\si^{(2)}(t,X_t^{(1)})\big\|^2\d t\bigg)\sup_{t\in [0,T]} h\big(|X_t^{(1)}| + |x^{(1)}-y^{(1)}|\big)\bigg\}  <\infty,$$ then for any strictly positive $f\in \B_b(\R^{m+d})$,
\beg{equation*}\beg{split}   P_{2T}\log f(y)\le & \log P_{2T} f(x)
+  \ff{K_T|x^{(1)}-y^{(1)}|^2}{\ll_T^2(1-\e^{-2K_TT})}\\
& +      \ff{\theta_T \e^{2\Theta_TT} \Psi_T(x^{(1)}, y^{(1)})}{2 }
    \bigg\{   |x^{(2)}-y^{(2)}|^2 +\ff{1-\e^{-2\Theta_TT}}{\Theta_T} \vp_T(|x^{(1)}-y^{(1)}|^2)\bigg\}. \end{split}\end{equation*}
\end{thm}

\

Because of Theorem \ref{T1.2}, we are now able to present  the log-Harnack inequality for the Gruschin semigroup on $\R^{m+d}$ with any $l>0$. Of course, one may also construct more general examples to illustrate  Theorem \ref{T1.2}.

\beg{cor}[Gruschin Semigroup] \label{C1.3}  Let $b^{(1)}=0, b^{(2)}=0, \si^{(1)}=I_{m\times m}$ and $\si^{(2)}(x^{(1)}) =|x^{(1)}|^lI_{d\times d}$ for some constant $l>0.$ Then there exists a constant $c>0$ such that
\beg{equation*}\beg{split} &P_{2T}\log f(y)\le \log P_{2T}f(x) + \ff{|x^{(1)}-y^{(1)}|^2}{2T} \\
&+\ff c {T^{l+1}} \Big(|x^{(1)}|^{2(l-1)^+} + |y^{(1)}|^{2(l-1)^+}+ T^{(l-1)^+}\Big)\Big(|x^{(2)}-y^{(2)}|^2 + 2T |x^{(1)}-y^{(1)}|^{2(l\land 1)}\Big)\end{split}\end{equation*} holds for all $T>0$ and $x,y\in \R^{m+d}.$
\end{cor}

In the next two sections, we will present proofs of Theorem \ref{T1.1}, Theorem \ref{T1.2} and Corollary \ref{C1.3} respectively.  In  the proof of Theorem \ref{T1.1}, the additional drifts constructed in the coupling are adapted so that the usual argument applies. However, in the proof of Theorem \ref{T1.2} the drift constructed for the coupling of the second component process is merely conditional adapted given $B^{(1)}$, a new trick is then introduced to derive the log-Harnack inequality.

\section{Proof of Theorem \ref{T1.1}}

Let $x=(x^{(1)}, x^{(2)}), y=(y^{(1)}, y^{(2)})$ and $T>0$ be fixed. The idea to establish a Harnack type inequality of $P_{2T}$ using  a  coupling by change of measure is as follows: construct two processes $X_t, Y_t$ and a probability density function $R$ such that   $X_{2T}=Y_{2T},  X_0=x, Y_0=y$,   and
$$P_{2T}f(x)= \E f(X_{2T}),\ \ P_{2T}f(y)= \E \big\{R f(Y_{2T})\big\},\ \ f\in\B_b(\R^{m+d}).$$ Then, by e.g. the Young inequality, for strictly positive $f$ one obtains
\beq\label{T} P_{2T}\log f(y) =\E \big\{R\log f(Y_{2T}) \big\}=  \E \big\{R\log f(X_{2T}) \big\}\le \E(R\log R) + \log P_{2T} f(x).\end{equation} This implies the log-Harnack inequality provided $\E (R\log R)<\infty.$

When the SDE is driven by an additive noise, this idea can be easily realized by adding a proper drift to the equation and using the Girsanov theorem. In the non-degenerate multiplicative noise case, the argument has been well modified in \cite{Wan3} by constructing a coupling with  singular additional drifts. For the present model, as the SDE is driven by a multiplicative noise with  a possibly  degenerate and singular coefficient, it is hard to follow the known ideas to construct a coupling in one go. What we will do in this paper is to  construct a  coupling in two steps, where the second step will be realized under the regular conditional probability given $B^{(1)}$: \beg{enumerate}
\item[(1)] We first construct a coupling $(X^{(1)}_t, Y_t^{(1)})$ by change of measure for  the first component of the process  such that $X_t^{(1)}=Y_t^{(1)}$ for $t\ge T$. This part is now standard as the first equation in (\ref{EX}) is driven by the  non-degenerate additive noise $\si^{(1)}(t)\d B_t^{(1)}.$
\item[(2)] Once $X_t^{(1)}=Y_t^{(1)}$ holds for $t\ge T$,  the equations for $X_t^{(2)}$ and $Y_t^{(2)}$ will have same noise part for $t\ge T$, so that we are able to construct a coupling by change of measure for them such that $X_{2T}^{(2)}=Y_{2T}^{(2)}.$
\end{enumerate}

\subsection{Construction of the coupling}

Throughout this section, we assume that {\bf (A1)}-{\bf (A3)} and condition (\ref{PS}) hold.
We first construct the Brownian motion $B_t$  as the coordinate process on the Wiener space
$(\OO,\F,\P)$, where  $$\OO= C([0,\infty); \R^{m+d})= C([0,\infty);\R^m)\times C([0,\infty);\R^d),$$
$\F$ is the Borel $\si$-field, $\P$ is the Wiener measure (i.e. the distribution of the $(m+d)$-dimensional Brownian motion starting at $0$). Let
$$B_t(\oo) = (B_t^{(1)}(\oo), B_t^{(2)}(\oo))= (\oo^{(1)}_t, \oo^{(2)}_t),\ \ \oo=(\oo^{(1)},\oo^{(2)})\in\OO, t\ge 0.$$ Then $B_t$ is the $(m+d)$-dimensional Brownian motion w.r.t. the natural filtration $(\F_t)_{t\ge 0}.$ Moreover, let $\F^{(1)}=\si(B_t^{(1)}: t\ge 0)$ and $\F^{(2)}_t=\si(B_s^{(2)}:0\le s\le t), t\ge 0.$ It is well known that the conditional regular probability $\P(\cdot|\F^{(1)})$ given $\F^{(1)}$ exists. This structure will enable us to first construct a coupling
$(X_t^{(1)}, Y^{(1)}_t)$ for the first component process up to time $T$ under probability $\P$, then construct a coupling $(X_t^{(2)}, Y_t^{(2)})$ for the second component process from time $T$ on under the regular conditional probability $\P(\cdot|\F^{(1)}).$ For any probability measure $\tilde{\mathbb{P}}$ on $(\OO,\F)$, we denote
by $\E_{\tilde{\mathbb{P}}}$ the expectation  w.r.t. $\tilde{\mathbb{P}}$. When $\tt \P=\P$, we simply denote the expectation by $\E$ as usual.

Let $X_t= (X_t^{(1)}, X_t^{(2)})$ solve the equation (\ref{EX}) with $X_0=x=(x^{(1)}, x^{(2)}).$ Given $Y_0=y=(y_1^{(1)},y_2^{(2)})\in\R^{m+d}$, we are going to construct $Y_t^{(1)}$ on $\R^m$ and $Y_t^{(2)}$ on $\R^d$ respectively, such that $Y_t^{(1)}=X_t^{(1)}$ for $t\ge T$ and $Y_{2T}^{(2)}= X_{2T}^{(2)}.$

\subsubsection{Construction of $Y_t^{(1)}$} Consider the equation
\beq\label{EY1} \d Y_t^{(1)}= b^{(1)}(t,Y_t^{(1)})\d t +\si^{(1)}(t) \d B_t^{(1)} - v_t^{(1)}\d t,\ \ Y_0^{(1)}=y^{(1)},\end{equation}
where
$$v_t^{(1)}:=\ff{2K_T|x^{(1)}-y^{(1)}|\e^{-K_Tt}(Y_t^{(1)}-X_t^{(1)})}{(1-\e^{-2K_TT})|X_t^{(1)}-Y_t^{(1)}|} \, 1_{\{X_t^{(1)}\ne Y_t^{(1)}\}},\ \ t\ge 0.$$
Obviously, the equation has a unique strong solution before the coupling time
$$\tau_1:= \inf\big\{t\ge 0: X_t^{(1)}=Y_t^{(1)}\big\}.$$ Then, letting $Y_t^{(1)}=X_t^{(1)}$ for $t\ge \tau_1$, we see that $(Y_t^{(1)})_{t\ge 0}$ is a strong solution to (\ref{EY1}). So, we can reformulate $v_t^{(1)}$ as
\beq\label{V1} v_t^{(1)}=\ff{2K_T|x^{(1)}-y^{(1)}|\e^{-K_Tt}(Y_t^{(1)}-X_t^{(1)})}{(1-\e^{-2K_TT})|X_t^{(1)}-Y_t^{(1)}|} \, 1_{[0,\tau_1)}(t),\ \ t\ge 0.\end{equation}

\beg{prp}\label{P1} For any $t\ge 0$,
\beq\label{X1} |X_t^{(1)}-Y_t^{(1)}|\le \ff{\e^{-K_Tt}-\e^{-K_T(2T-t)}}{1-\e^{-2K_TT}} |x^{(1)}-y^{(1)}|1_{[0,T]}(t)\le |x^{(1)}-y^{(1)}|1_{[0,T]}(t). \end{equation} Consequently, $\tau_1\le T$ and $X_t^{(1)}=Y_t^{(1)}$ for $t\ge T$.\end{prp}

\beg{proof} By {\bf (A2)} and (\ref{V1}), we have
$$\d |X_t^{(1)}- Y_t^{(1)}| \le \bigg(K_T |X_t^{(1)} -Y_t^{(1)}|-\ff {2K_T |x^{(1)}- y^{(1)}|\e^{-K_Tt}}{1-\e^{-2K_TT}}\bigg)\d t,\ \ t\in [0,\tau_1)\cap [0,T].$$ Then
$$|X_t^{(1)}-Y_t^{(1)}|\le \ff{\e^{-K_Tt}-\e^{-K_T(2T-t)}}{1-\e^{-2K_TT}} |x^{(1)}-y^{(1)}|,\ \ t\in [0,\tau_1)\land [0,T].$$
This implies that $\tau_1\le T$  and   also (\ref{X1}) since $X_t^{(1)}=Y_t^{(1)}$ for $t\ge \tau_1.$
\end{proof}

To formulate (\ref{EY1}) as the first equation in (\ref{EX}), we let
$$\tt B_t^{(1)}= B_t^{(1)} - \int_0^{t} \xi^{(1)}(s)\d s,\ \ \xi^{(1)}(s):=   \si^{(1)}(t)^{-1} v_t^{(1)},\ \ t \ge 0.$$  From {\bf (A1)} and
(\ref{V1})  we see that $\xi^{(1)}(s)$ is bounded and adapted. So,  by the Girsanov theorem, $\tt B_t$ is an $m$-dimensional Brownian motion under the probability measure $\Q^{(1)}:= R_1(T)\P$, where
$$R_1(t):= \exp\bigg[\int_0^t \<\xi^{(1)}(s), \d B_s^{(1)}\> -\ff 1 2 \int_0^t |\xi^{(1)}(s)|^2\d s\bigg],\ \ t\ge 0$$ is a martingale.
Obviously,  (\ref{EY1}) can be formulated as
\beq\label{EY1'} \d Y_t^{(1)}= b^{(1)}(t,Y_t^{(1)})\d t + \si^{(1)}(t) \d\tt B_t^{(1)},\ \ Y_0^{(1)}=y^{(1)}.\end{equation}
As shown in (\ref{T}), for the log-Harnack inequality we need to estimate the entropy of $R_1:=R_1(T)$.

\beg{prp}\label{P2} Let $R_1=R_1(T).$ Then
\beq\label{ENT1} \E\big\{R_1\log R_1\big\}\le \ff {K_T|x^{(1)}-y^{(1)}|^2}{\ll_T^2(1-\e^{-2K_TT})}.\end{equation}
\end{prp}

\beg{proof} By $\tau_1\le T$, {\bf (A1)}, (\ref{V1}), we have
\beq\label{D1} \int_0^T |\si^{(1)}(t)^{-1} v_t^{(1)}|^2\d t\le  \ff {2K_T|x^{(1)}-y^{(1)}|^2}{\ll_T^2(1-\e^{-2K_TT})}.\end{equation}
Then, it follows from (\ref{EY1})    and the definition of $R_1$ that
\beg{equation*}\beg{split} &\E\big\{R_1\log R_1\big\}= \E_{\Q^{(1)}}\log R_1 \\
&= \ff 1 2 \E_{\Q^{(1)}} \int_0^{T} |\si^{(1)}(t)^{-1} v_t^{(1)}|^2\d t \le \ff {K_T|x^{(1)}-y^{(1)}|^2}{\ll_T^2(1-\e^{-2K_TT})}.\end{split}\end{equation*}
\end{proof}

\subsubsection{Construction of $Y_t^{(2)}$}

Consider the equation
\beq\label{EY2} \d Y_t^{(2)}=  b^{(2)} (t,Y_t) \d t + \si^{(2)}(t,Y_t^{(1)}) \d B_t - v_t^{(2)}\d t,\ \ Y_0^{(2)}=y^{(2)},\end{equation}
where
$$v_t^{(2)}:= \ff{2\Theta_{2T}|X_T^{(2)}-Y_T^{(2)}|\e^{-\Theta_{2T}t}(Y_t^{(2)}-X_t^{(2)})}{(\e^{-2\Theta_{2T}T}-\e^{-4\Theta_{2T}T})|X_t^{(2)}-Y_t^{(2)}|}1_{\{t\ge T, X_t^{(2)}\ne Y_t^{(2)}\}},\ \ t\ge 0.$$
As $Y_t^{(1)}$ is now fixed, it is easy to see that (\ref{EY2}) has a unique solution before time
$$\tau_2:=\inf\{t\ge T: X_t^{(2)}= Y_t^{(2)}\}.$$ Letting $Y_t^{(2)}=X_t^{(2)}$ for $t\ge \tau_2$, we see that $(Y_t^{(2)})_{t\ge 0}$ solves the equation (\ref{EY2}). Thus,
\beq\label{V2}  v_t^{(2)}= \ff{2\Theta_{2T}|X_T^{(2)}-Y_T^{(2)}|\e^{-\Theta_{2T}t}(Y_t^{(2)}-X_t^{(2)})}{(\e^{-2\Theta_{2T}T}-\e^{-4\Theta_{2T}T})|X_t^{(2)}-Y_t^{(2)}|}1_{[T,\tau_2)}(t),\ \ t\ge 0.\end{equation}

\beg{prp}\label{P3} For any $t\ge T$,
\beq\label{X2} |X_t^{(2)}-Y_t^{(2)}|\le \ff{\e^{-\Theta_{2T}(t-T)}-\e^{-\Theta_{2T}(3T-t)}}{1-\e^{-2\Theta_{2T}T}}|X_T^{(2)}-Y_T^{(2)}|1_{[T,2T]}(t).\end{equation} \end{prp}

\beg{proof}  Since $\vp_\cdot(0)=0$ and $X_t^{(1)}=Y_t^{(1)}$ for $t\ge T$, by {\bf (A3)}, (\ref{V2})  and It\^o's formula we obtain
$$\d |X_t^{(2)}-Y_t^{(2)}|\le \bigg(\Theta_{2T}|X_t^{(2)}-Y_t^{(2)}|-\ff{2\Theta_{2T}|X_T^{(2)}-Y_T^{(2)}|\e^{-\Theta_{2T}t}}{\e^{-2\Theta_{2T}}- \e^{-4\Theta_{2T}T}}\bigg)\d t,\ \  t\in [T,\tau_2)\cap [T,2T].$$ This implies (\ref{X2}) for $t\in [T,\tau_2)\cap [T,2T].$ Therefore, $\tau_2\le 2T$ and (\ref{X2}) holds for all $t\ge T$.
\end{proof}

To formulate (\ref{EY2}) as the second equation in (\ref{EX}), we need to make use of the Girsanov theorem to get rid of the additional drift. To this end, let
$$\xi^{(2)}(s)= \si^{(2)}(s,Y_s^{(1)})^{-1} v_s^{(2)},\ \ \ s\in [T,2T]$$ and
$$R_2(t)= \exp\bigg[\int_T^t \<\xi^{(2)}(s), \d B_s^{(2)}\>-\ff 1 2 \int_T^t |\xi^{(2)}(s)|^2\d s\bigg],\ \ t\in [T,2T].$$
Since $B^{(2)}_t$ is independent of $\F^{(1)}$, the following result ensures that $\{R_2(t)\}_{t\in [T,2T]}$ is a uniformly integrable $\F^{(2)}_t$-martingale under $\P(\cdot|\F^{(1)}).$

\beg{prp}\label{P4} Under $\P(\cdot|\F^{(1)}), \{R_2(t)\}_{t\in [T,2T]}$ is an $\F_t^{(2)}$-martingale   and $R_2:=R_2(2T)$ satisfies
\beq\label{W*0} \beg{split}  &\E_{\P(\cdot|\F^{(1)})} \big\{R_2\log R_2\big\}
 \le  \bigg(\int_T^{2T} \ff{2\Theta_{2T}^2\e^{-2\Theta_{2T}t}\|\si^{(2)}(t,Y_t^{(1)})^{-1}\|^2 }{(\e^{-2\Theta_{2T}T}-\e^{-4\Theta_{2T}T})^2}\d t\bigg)\\
 &\times \bigg(\e^{2\Theta_TT}|x^{(2)}-y^{(2)}|^2 +\ff{\e^{2\Theta_TT}-1}{\Theta_T}   \vp_T(|x^{(1)}-y^{(1)}|^2)\bigg)\sup_{t\in [0,T]} h(|Y_t^{(1)}|+|x^{(1)}-y^{(1)}|).
 \end{split} \end{equation}\end{prp}

\beg{proof} We make use of an approximation argument. Let $\xi_n^{(2)}(s)= \xi^{(2)}(s)1_{\{|\xi^{(2)}(s)|\le n\}},$  and let
$$R_{2,n}(t)= \exp\bigg[\int_T^t \<\xi^{(2)}_n(s), \d B_s^{(2)}\>-\ff 1 2 \int_T^t |\xi^{(2)}_n(s)|^2\d s\bigg],\ \ n\ge 1, t\in [T,2T].$$ Then
$\{R_{2,n}(t)\}_{t\in [T,2T]}$ is an $\F_t^{(2)}$-martingale under $\P(\cdot|\F^{(1)})$.   So, it remains to show that
\beq\label{W*1} \beg{split} &\E_{\P(\cdot|\F^{(1)})} \big\{R_{2,n}\log R_{2,n}\big\}  (t)
 \le  \bigg(\int_T^{2T} \ff{2\Theta_{2T}^2\e^{-2\Theta_{2T}t}\|\si^{(2)}(t,Y_t^{(1)})^{-1}\|^2 }{(\e^{-2\Theta_{2T}T}-\e^{-4\Theta_{2T}T})^2}\d t\bigg)\\
 &\times \bigg(\e^{2\Theta_TT}|x^{(2)}-y^{(2)}|^2 +\ff{\e^{2\Theta_TT}-1}{\Theta_T}   \vp_T(|x^{(1)}-y^{(1)}|^2)\sup_{t\in [0,T]} h(|Y_t^{(1)}|+|x^{(1)}-y^{(1)}|) \bigg)
  \end{split} \end{equation}holds for all $t\in [T,2T]$ and $n\ge 1.$
Let $\Q_{2,n} = R_{2,n}(2T) \P(\cdot|\F^{(1)}).$ By the Girsanov theorem, under $\Q_{2,n}$ the process
$$\tt B^{(2)}_t:= B_t^{(2)} -\int_T^{T\lor t} \xi_n^{(2)}(s)\d s,\ \ \ t\in [0,2T]$$ is a $d$-dimensional Brownian motion.
Then, by the definition of $\xi_n^{(2)}(s)$ and (\ref{V2}), we have
\beg{equation}\label{W*3} \beg{split} &\E_{\P(\cdot|\F^{(1)})} \big\{R_{2,n}\log R_{2,n}\big\}  (2T) = \E_{\Q_{2,n}} \log R_{2,n}(2T) = \ff 1 2 \int_T^{2T} \E_{\Q_{2,n}}|\xi_n^{(2)}(s)|^2\d s\\
&\le     \bigg( \int_T^{2T} \ff{2\Theta_{2T}^2\e^{-2\Theta_{2T}t}\|\si^{(2)}(t,Y_t^{(1)})^{-1}\|^2}{(\e^{-2\Theta_{2T}T}-\e^{-4\Theta_{2T}T})^2}\,\d t\bigg)\E_{\P(\cdot|\F^{(1)})}\big\{R_{2,n} (2T)|X_T^{(2)}-Y_T^{(2)}|^2\big\}.\end{split}\end{equation}
Since $\{R_{2,n}(t)\}_{t\in [T,2T]}$ is an $\F_t^{(2)}$-martingale under $P(\cdot|\F^{(1)})$, and $R_{2,n}(T)=1$,  we have
\beq\label{W*4} \E_{\P(\cdot|\F^{(1)})} \big\{R_{2,n} (2T) |X_T^{(2)}-Y_T^{(2)}|^2\big\}= \E_{\P(\cdot|\F^{(1)})}  |X_T^{(2)}-Y_T^{(2)}|^2.\end{equation}
Finally, by {\bf (A3)}, (\ref{X1}) and It\^o's formula, we obtain
\beg{equation*}\beg{split} \d |X_t^{(2)}-Y_t^{(2)}|^2 \le &2 \big\{\Theta_T|X_t^{(2)}-Y_t^{(2)}|^2 +\vp_T(|x^{(1)}-y^{(1)}|^2)h(|Y_t^{(1)}|+|x^{(1)}-y^{(1)}|) \big\}\d t\\
& + 2 \big\<X_t^{(2)}-Y_t^{(2)}, \{\si^{(2)}(t,X_t^{(1)})-\si^{(2)}(t,Y_t^{(1)})\}\d B_t^{(2)}\big\>,\ \ t\in [T,2T].\end{split}\end{equation*}   Since $h\ge 1$, this implies
\beq\label{D2} \beg{split} &\E_{\P(\cdot|\F^{(1)})}  |X_T^{(2)}-Y_T^{(2)}|^2\\
&\le \bigg(\e^{2\Theta_TT}|x^{(2)}-y^{(2)}|^2 +\ff{\e^{2\Theta_TT}-1}{\Theta_T} \vp_T(|x^{(1)}-y^{(1)}|^2)\bigg)\sup_{t\in [0,T]} h(|Y_t^{(1)}|+|x^{(1)}-y^{(1)}|).\end{split}\end{equation}
Combining this with (\ref{W*3}) and (\ref{W*4}), we prove (\ref{W*1}).
\end{proof}

\beg{proof}[Proof of Theorem \ref{T1.1}]  Let $X_t=(X_t^{(1)}, X_t^{(2)})$ and $Y_t= (Y_t^{(1)}, Y_t^{(2)})$ be constructed above. Let $R=  R_1 R_2$.
By Propositions \ref{P1}, \ref{P2}, \ref{P3} and \ref{P4}, we have $X_{2T}=Y_{2T},$
  $\E_{\P(\cdot|\F^{(1)})}R_2  =1,$  and  noting that the distribution of $Y^{(1)}$ under $R_1\P$ coincides with that of $X^{(1)}(y^{(1)})$ under $\P$,
\beg{equation*}\beg{split} &\E\{R\log R\}= \E\Big\{(R_1\log R_1) \E_{\P(\cdot|\F^{(1)})} R_2 \Big\}+\E \Big\{R_1  \E_{\P(\cdot|\F^{(1)})} (R_2 \log R_2) \Big\}  \\
&\le  \ff{K_T|x^{(1)}-y^{(1)}|^2}{\ll_T^2(1-\e^{-2K_TT})}
 +    \E^{y^{(1)}}\bigg\{\bigg(\int_T^{2T} \ff{2\Theta_{2T}^2\e^{-2\Theta_{2T}t}\|\si^{(2)}(t,X_t^{(1)})^{-1}\|^2 }{(\e^{-2\Theta_{2T}T}-\e^{-4\Theta_{2T}T})^2}\d t\bigg)\\
 &\quad \times \bigg( \e^{2\Theta_TT} |x^{(2)}-y^{(2)}|^2 +\ff{\e^{2\Theta_TT}-1}{\Theta_T} \vp_T(|x^{(1)}-y^{(1)}|^2)\bigg)\sup_{t\in [0,T]} h(|X_t^{(1)}|+|x^{(1)}-y^{(1)}|) \bigg\}\\
 &\le \ff{K_T|x^{(1)}-y^{(1)}|^2}{\ll_T^2(1-\e^{-2K_TT})}
   +  \ff{\Theta_{2T}\e^{2\Theta_TT}\psi_T(x^{(1)}, y^{(1)})}{\e^{-2\Theta_{2T}T}-\e^{-4\Theta_{2T}T}}
  \bigg\{   |x^{(2)}-y^{(2)}|^2 +\ff{1-\e^{-2\Theta_TT}}{\Theta_T} \vp_T(|x^{(1)}-y^{(1)}|^2)\bigg\}.
\end{split}\end{equation*}
 Therefore, the desired log-Harnack inequality   follows from (\ref{T}), since
  under the probability measure $\Q:= R\P$
 $$\tt B_t:= B_t +\int_0^t (\xi^{(1)}(s), \xi^{(2)}(s))\d s,\ \ t\ge 0$$ is a Brownian motion on $\R^{m+d}$, and $Y_t$ with $Y_0=y$ solves the equation
 $$\beg{cases} \d Y_t^{(1)}= b^{(1)}(t, Y_t^{(1)}) \d t + \si^{(1)}(t)\, \d \tt B_t^{(1)},\\
\d \tt Y_t^{(2)}= b^{(2)} (t,Y_t)\d t + \si^{(2)}(t,Y_t^{(1)})\,\d  \tt B_t^{(2)},\end{cases}$$ so that $P_{2T}f(y)= \E_\Q f(Y_{2T})= \E \big\{R f(Y_{2T})\big\}
=\E \big\{R f(X_{2T})\big\}.$
  \end{proof}

\section{Proofs of Theorem \ref{T1.2} and Corollary \ref{C1.3}}

\beg{proof}[Proof of Theorem \ref{T1.2}] Let $X_t= (X_t^{(1)}, X_t^{(2)})$  and $Y_t^{(1)}$ be constructed   in the last section. We now modify the construction of $Y_t^{(2)}$ in terms of the condition $\Psi_T<\infty$.
 Let
 $$\eta_t =\si^{(2)}(t, Y_t^{(1)})^*\e^{A^*(T-t)}  Q_T^{-1} (Y_T^{(2)}-X_T^{(2)})1_{[T,2T]}(t),\ \ t\ge 0.$$ Let $Y_t^{(2)}$ solve the equation
 \beq\label{D1} \d Y_t^{(2)}= b^{(2)}(t, Y_t) \d t +\si^{(2)}(t, Y_t^{(1)}) \big\{\d B_t^{(2)} - \eta_t\d t\big\},\ \ Y_0^{(2)}=y^{(2)}. \end{equation}
 Since under $\P(\cdot|\F^{(1)})$ the processes $X_t^{(1)}$ and $Y_t^{(1)}$ are fixed and $B_t^{(2)}$ is a $d$-dimensional Brownian motion,
 by {\bf (A3)}
 this equation has a unique solution. Since $X_t^{(1)}= Y_t^{(1)}$ for $t\ge T,$  for the present $b^{(2)}$ we have
 $b^{(2)}(t,X_t)- b^{(2)}(t, Y_t)= A(X_t^{(2)}- Y_t^{(2)})$ for $t\ge T$.  So,
 $$ X_{2T}^{(2)}-Y_{2T}^{(2)} =  \e^{AT} (X_T^{(2)}-Y_T^{(2)}) +\int_T^{2T} \e^{A(2T-t)} \si^{(2)}(t,Y_t^{(1)})\eta_t \d t=0$$ as $Y_t^{(1)}=X_t^{(1)}$ for $t\ge T$.
 Therefore, $X_{2T}=Y_{2T}.$  Moreover, let
 $$\tt R_2 = \exp\bigg[\int_T^{2T} \<\eta_t, \d B_t^{(2)}\> -\ff 1 2 \int_T^{2T} |\eta_t|^2\d t\bigg].$$ Following the proof of Proposition \ref{P4} and using (\ref{D2}), we obtain
 \beg{equation*}\beg{split} &\E_{\P(\cdot|\F^{(1)})} \big\{\tt R_2\log\tt R_2\big\} = \ff 1 2 \int_T^{2T} \E_{\tt R_2 \P(\cdot|\F^{(1)})} |\eta_t|^2\d t\\
 &\le \ff{\theta_T} 2 \Big(\E_{\P(\cdot|\F^{(1)})} \|X_T^{(2)}-Y_T^{(2)}\|^2 \Big)   \|Q_T^{-1}\|^2 \int_T^{2T} \|\si^{(2)}(t, Y_t^{(1)}) \|^2 \d t\\
 &\le \ff{\theta_T} 2 \bigg\{\e^{2\Theta_TT}|x^{(2)}-y^{(2)}|^2 +\ff{\e^{2\Theta_TT}-1}{\Theta_T} \vp_T(|x^{(1)}-y^{(1)}|^2)\bigg\}\sup_{t\in [0,T]} h(|Y_t^{(1)}|+|x^{(1)}-y^{(1)}|)\\
 &\qquad\times        \|Q_T^{-1}\|^2 \int_T^{2T} \|\si^{(2)}(t, Y_t^{(1)}) \|^2 \d t.
 \end{split}\end{equation*}
Repeating the proof of Theorem \ref{T1.1}  and  using this inequality to replace (\ref{W*0}), we obtain
\beq\label{AB0} \beg{split}&\E\big\{(R_1\tt R_2)\log (R_1\tt R_2)\big\} \le \ff{K_T|x^{(1)}-y^{(1)}|^2}{\ll_T^2(1-\e^{-2K_TT})}\\
& +      \ff{\theta_T \e^{2\Theta_TT} \Psi_T(x^{(1)}, y^{(1)})}{2 }
    \bigg\{   |x^{(2)}-y^{(2)}|^2 +\ff{1-\e^{-2\Theta_TT}}{\Theta_T} \vp_T(|x^{(1)}-y^{(1)}|^2)\bigg\}.
\end{split} \end{equation}
Since $B_t^{(2)}$  is a $d$-dimensional Brownian motion under $\P(\cdot|\F^{(1)})$,   by the Girsanov theorem, under $\tt R_2\P(\cdot|\F^{(1)})$ the process
$$\tt B_t^{(2)}:= B_t^{(2)}- \int_T^t\eta_s\d s,\ \ t\in [T,2T]$$ is a $d$-dimensional Brownian motion. Noting that
$$Y_t^{(2)}= Y_T^{(2)} +\int_T^t b^{(2)}(s, Y_s)\d s + \int_T^t \si^{(2)}(s, Y_s^{(1)}) \d \tt B_s^{(2)},\ \ t\in [T,2T],$$
we see that the distribution of $Y_{2T}^{(2)}$ under $\tt R_2\P(\cdot|\F^{(1)})$ coincides with that of $\tt Y_{2T}^{(2)}$ under $\P(\cdot|\F^{(1)})$, where
$$ \tt Y_t^{(2)} =\beg{cases} Y_t^{(2)}, &\text{if}\ t\in [0,T],\\
Y_T^{(2)} +\int_T^t b^{(2)}(s,Y_s)\d s + \int_T^t \si^{(2)}(s, Y_s^{(1)})\d B_s^{(2)}, &\text{if}\ t\in [T, 2T].\end{cases}$$
Therefore,
$$\E_{\P(\cdot|\F^{(1)})} \big\{\tt R_2 \log f(Y_{2T})\big\}= \E_{\P(\cdot|\F^{(1)})} \big\{\log  f(Y_{2T}^{(1)}, \tt Y_{2T}^{(2)})\big\}.$$
Combining this with $X_{2T}=Y_{2T}$, we obtain
\beg{equation}\label{AB}\beg{split} &\E\big\{R_1\tt R_2\log f(X_{2T})\big\}= \E\big\{R_1\tt R_2\log f(Y_{2T})\big\} =\E\Big(R_1 \E_{\P(\cdot|\F^{(1)})} \big\{\tt R_2 \log f(Y_{2T})\big\}\Big)\\
&= \E\Big(R_1 \E_{\P(\cdot|\F^{(1)})} \big\{\log  f(Y_{2T}^{(1)}, \tt Y_{2T}^{(2)})\big\}\Big) = \E\big\{R_1\log f(Y_{2T}^{(1)}, \tt Y_{2T}^{(2)})\big\}.\end{split}\end{equation}
Moreover, again by the Girsanov theorem, under $R_1\P$ the process $(\tt B_t^{(1)}, B_t^{(2)})_{t\in [0,2T]}$ is a $(d+m)$-dimensional Brownian motion, recall that
$$\tt B_t^{(1)}= B_t^{(1)}-\int_0^{T\land t} \xi^{(1)}(s)\d s,\ \ t\in [0,2T].$$ Noting that $(Y_t^{(1)},\tt Y_t^{(2)})$ solves the equation
$$\beg{cases} \d Y_t^{(1)}= b^{(1)}(t, Y_t^{(1)}) \d t + \si^{(1)}(t)\, \d \tt B_t^{(1)},\ Y_0^{(1)}=y^{(1)},\\
\d \tt Y_t^{(2)}= b^{(2)} (t,Y_t^{(1)}, \tt Y_t^{(2)})\d t + \si^{(2)}(t,Y_t^{(1)})\,\d  B_t^{(2)},\  \tt Y_0^{(2)}= y^{(2)},\end{cases}$$
we conclude that the distribution of $(Y_{2T}^{(1)}, \tt Y_{2T}^{(2)})$ under $R_1\P$ coincides with that of $X_{2T}(y)$ under $\P$. Therefore, it follows from (\ref{AB}) and the Young inequality that
\beg{equation*}\beg{split} P_{2T}\log f(y)& = \E\big\{R_1\log f(Y_{2T}^{(1)}, \tt Y_{2T}^{(2)})\big\} = \E\big\{R_1\tt R_2\log f(X_{2T})\big\}\\
&\le \log P_{2T}f(x) +\E\{(R_1\tt R_2)\log (R_1\tt R_2)\big\}.\end{split}\end{equation*}
Combining this with (\ref{AB0}) we complete the proof.
\end{proof}

\beg{proof}[Proof of Corollary \ref{C1.3}] It is easy to see that {\bf (A1)}-{\bf (A3)} hold for $\ll=1, K=\Theta=0, \varphi (r)=r^{l\land 1}$ and $h(r)= c_1 \lor r^{2(l-1)^+}$ for some constant $c_1\ge 1.$  Moreover,
$$Q_T=I_{d\times d} \int_T^{2T} |B_t^{(1)}+ x^{(1)}|^{2l}\d t$$ is invertible and
$$\|Q_T^{-1}\|^2 \int_T^{2T} \|\si^{(2)}(X^{(1)})\|^2\d t= \ff 1 {\int_T^{2T} |B_t^{(1)} +x^{(1)}|^{2l}\d t}.$$ Then, using  the fact that for any $r\ge 0$, $$\E \sup_{t\in [0,T] } |B_t^{(1)}+x^{(1)}|^{2r}\le c(r) (|x^{(1)}|^{2r} + T^r)$$ holds for some constant $c(r)>0,$  and noting that \cite[Lemma 3.1]{W12} implies
\beg{equation*}\beg{split} &\E\bigg( \int_T^{2T} |B_t^{(1)}+x^{(1)}|^{2l}\d t \bigg)^{-2}\\
&= \E\bigg\{ \E\bigg(\bigg( \int_0^{T} \big|(B_{T+t}^{(1)}-B^{(1)}_T)+(B_T^{(1)}+x^{(1)}\big)\big|^{2l}\d t \bigg)^{-2}\bigg|B_T^{(1)}\bigg)\bigg\} \le \ff {C} {T^{2(l+1)}}\end{split}\end{equation*} for some constant $C>0,$
we conclude that
\beg{equation*}\beg{split} \Psi_T(x^{(1)}, y^{(1)}) &\le   \bigg( \E\sup_{t\in [0,T] }h(|B_t^{(1)}+x^{(1)}|+ |x^{(1)}-y^{(1)}|)^2\bigg)^{\ff 1 2}
\bigg(\E\bigg( \int_T^{2T} |B_t^{(1)}+x^{(1)}|^{2l}\d t \bigg)^{-2}\bigg)^{\ff 1 2} \\
&\le \ff c {T^{l+1}} \Big( |x^{(1)}|^{2(l-1)^+} + |y^{(1)}|^{2(l-1)^+}+ T^{(l-1)^+}\Big)\end{split}\end{equation*} holds for some constant $c>0.$ Therefore, the desired
log-Harnack inequality follows from Theorem \ref{T1.2}.
\end{proof}

  \beg{thebibliography}{99}

 \bibitem{ATW06} M. Arnaudon, A. Thalmaier, F.-Y. Wang,
  \emph{Harnack inequality and heat kernel estimates
  on manifolds with curvature unbounded below,} Bull. Sci. Math. 130(2006), 223--233.

\bibitem{ATW09} M. Arnaudon, A. Thalmaier, F.-Y. Wang,
  \emph{Gradient estimates and Harnack inequalities on non-compact Riemannian manifolds,}
   Stoch. Proc. Appl. 119(2009), 3653--3670.

   \bibitem{BBBC} D. Bakry,  F. Baudoin,  M. Bonnefont, D. Chafa\"{\i}, \emph{ On gradient bounds for the heat kernel on the Heisenberg group,} J. Funct. Anal. 255 (                    2008),   1905--1938.

\bibitem{BB1} F. Baudoin, M. Bonnefont, \emph{Log-Sobolev inequalities for subelliptic operators satisfying a generalized curvature dimension inequality,} J. Funct. Anal. 262(2012), 2646--2676.


\bibitem{BB3} F. Baudoin, N. Garofalo, \emph{Curvature-dimension inequalities and Ricci lower bounds for sub-Riemannian manifolds with transverse symmetries,} arXiv:1101.3590.

\bibitem{BGM} F. Baudoin, M. Gordina, T. Melcher, \emph{Quasi-invariance for heat kernel measures on sub-Riemannian infinite-dimensional Heisenberg groups,} arXiv:1108.1527.

\bibitem{B} J. M. Bismut, \emph{Large Deviations and the
Malliavin Calculus,} Boston: Birkh\"auser, MA, 1984.

\bibitem{Hand} A. N. Borodin, P. Salminen, \emph{Handbook of Brownian Motion - Facts and Formulae,} Birkh\"auser, Berlin, 1996.

\bibitem{DM} B. K. Driver, T. Melcker, \emph{Hypoelliptic heat kernel inequalities on the Heisenberg group,} J. Funct. Anal. 221(2005), 340--365.

\bibitem{GW} A. Guillin, F.-Y. Wang,
  \emph{Degenerate Fokker-Planck equations : Bismut formula, gradient estimate  and Harnack inequality,}  J. Diff. Equat. 253(2012), 20--40

\bibitem{Li} H.-Q. Li, \emph{Estimation optimale du gradient du semi-groupe de la chaleur sur le groupe de Heisenberg,} J. Funct. Anal. 236(2006), 369--394.

\bibitem{L} W. Liu, \emph{Harnack inequality and applications for stochastic evolution equations with monotone drifts,}  J. Evol. Equ. 9 (2009),  747--770.

\bibitem{LW} W.  Liu, F.-Y.  Wang,  \emph{Harnack inequality and strong Feller property for stochastic fast-diffusion equations,} J. Math. Anal. Appl. 342(2008), 651--662.

\bibitem{O} S.-X. Ouyang, \emph{Harnack inequalities and applications for multivalued stochastic evolution equations,}  Inf. Dimen.
Anal. Quant. Probab. Relat. Topics.

\bibitem{ORW} S.-X. Ouyang, M. R\"ockner, F.-Y. Wang, \emph{Harnack inequalities and applications for Ornstein-Uhlenbeck semigroups with jump,}
 Potential Anal. 36(2012), 301--315.

\bibitem{RW09} M. R\"ockner, F.-Y. Wang, \emph{Log-Harnack  inequality for stochastic differential equations in Hilbert spaces and its consequences, } Infin. Dimens. Anal. Quant. Probab.  Relat. Topics 13(2010), 27--37.

\bibitem{W97} F.-Y. Wang,  \emph{Logarithmic Sobolev inequalities on noncompact Riemannian manifolds,} Probab. Theory Relat. Fields 109(1997), 417-424.
\bibitem{W10} F.-Y. Wang \emph{Harnack inequalities on manifolds with boundary and applications,}    J.
Math. Pures Appl.   94(2010), 304--321.

\bibitem{Wan3} F.-Y. Wang,  \emph{Harnack inequality for SDE with multiplicative noise and extension to
               Neumann semigroup on nonconvex manifolds,} Ann. Probab. 39(2011), 1447-1467.

 \bibitem{W12a} F.-Y. Wang, \emph{Derivative formula and gradient estimates for  Gruschin type semigroups,} to appear in J. Theo. Probab.

\bibitem{W12} F.-Y. Wang, \emph{Generalized curvature condition for subelliptic diffusion processes,} arXiv: 1202.0778.

  \bibitem{WWX10}  F.-Y. Wang, J.-L. Wu and L. Xu, \emph{Log-Harnack inequality for stochastic Burgers equations and
applications},     J. Math. Anal. Appl. 384(2011), 151--159.

\bibitem{WY10} F.-Y. Wang, C. Yuan, \emph{Harnack inequalities for functional SDEs with multiplicative noise and applications,}   Stoch. Proc. Appl.    121(2011), 2692--2710.

\bibitem{WX10} F.-Y. Wang, L. Xu, \emph{Derivative formula and applications for hyperdissipative stochastic Navier-Stokes/Burgers equations,} to appear in Inf. Dim. Quant. Probab. Relat. Topics arXiv:1009.1464.

\bibitem{WZ} F.-Y. Wang, X.-C. Zhang, \emph{Derivative formula and applications  for degenerate diffusion semigroups,} arXiv1107.0096.

\bibitem{X} L. Xu, \emph{A modified log-Harnack inequality and asymptotically strong Feller property,}  J. Evol. Equ. 11 (2011),  925--942.

\bibitem{Z} X.-C. Zhang, \emph{Stochastic flows and Bismut formulas for stochastic Hamiltonian systems,} Stoch. Proc. Appl. 120(2010), 1929--1949.

\end{thebibliography}
\end{document}